\RequirePackage[leqno]{amsmath}
\documentclass[leqno]{amsart}
\usepackage{amsmath}
\usepackage[pagewise]{lineno}
\makeatletter
\usepackage[utf8]{inputenc}     
\usepackage[T1]{fontenc}  
\usepackage{multicol,multirow}
\usepackage{color, graphics}
\usepackage{pgfplots}
\usepackage[output-decimal-marker={,},exponent-product=\cdot]{siunitx}
\usepackage{graphicx, multicol, latexsym, amsmath, amssymb}
\usepackage{amsfonts}
\usepackage{enumitem}
\setlist[itemize]{topsep=0pt,after=\vspace{1.5\baselineskip}}
 \usepackage[colorlinks=true]{hyperref}
\usepackage{empheq}
\usepackage[T1]{fontenc}
\usepackage{tabularx}
\usepackage[leftcaption]{sidecap}
\sidecaptionvpos{figure}{m}
\usepackage{tikz}
\usepackage{subfigure}
\usepackage[margin=0.5in]{geometry}
\usetikzlibrary{backgrounds,automata}

\newcount\rc@count
\let\rc@clearconstantlist\empty
\newcommand\rc@clearconstant[1]{\global\expandafter\let\csname rc@const@#1\endcsname\undefined}
\newcommand\resetconstants[1]{%
    \def\rc@constname{#1}
    \global\rc@count=1\relax 
    \bgroup 
        \let\\\rc@clearconstant 
        \rc@clearconstantlist
        \global\let\rc@clearconstantlist\empty 
    \egroup
}
\newcommand\const[1]{%
    \@ifundefined{rc@const@#1}{%
        \expandafter\xdef\csname rc@const@#1\endcsname{%
           \noexpand\rc@useconst{\rc@constname}{\the\rc@count}%
        }%
        \g@addto@macro\rc@clearconstantlist{\\{\mathrm{#1}}}%
        \global\advance\rc@count1\relax
    }{}%
    \csname rc@const@#1\endcsname
}
\newcommand\rc@useconst[2]{{#1}\textsubscript{#2}}
\setlist[itemize]{noitemsep, topsep=0pt}

\def\R{\mathbb R} \def\N{\mathbb N}

\def\R{\mathbb R} \def\N{\mathbb N} 
\def\TM{T_{max}} 
\def
\@cite
#1#2{[{{\bfseries #1}\if@tempswa , #2\fi}]}
\makeatother

\newtheorem{theorem}{Theorem}[section]

\newtheorem{lemma}[theorem]{Lemma}
\newtheorem{proposition}{Proposition}

\newtheorem{remark}{Remark}

\newcounter{cnstcnt}

\title[Boundedness in an attraction-repulsion chemotaxis system with nonlinear productions] 
      {
Addendum to the paper ``Refined criteria toward boundedness in an attraction-repulsion chemotaxis system with nonlinear productions'' }
\author[S. Frassu and G. Viglialoro]{}
\keywords{Chemotaxis, Global existence, Boundedness, Nonlinear production. \\
	\textit{$^*$Corresponding author}: silvia.frassu@unica.it}
\begin{document}
\maketitle

\centerline{\scshape Silvia Frassu$^{*}$ \and Giuseppe Viglialoro}
\medskip
{
 \medskip
 \centerline{Dipartimento di Matematica e Informatica}
 \centerline{Universit\`{a} degli Studi di Cagliari}
 \centerline{Via Ospedale 72, 09124. Cagliari (Italy)}
 \medskip
}
\bigskip
\begin{abstract}
These notes aim to provide a deeper insight on the specifics of the paper ``Refined criteria toward boundedness in an attraction-repulsion chemotaxis system with nonlinear productions'' by A. Columbu, S. Frassu and G. Viglialoro [\textit{Appl. Anal.} 2024, 103:2, 415--431].

\end{abstract}
\resetconstants{c}
\section{Aim of the paper}
In this report we focus on \cite[Theorem 2.2]{CFVApplAnal2024} where
an attraction-repulsion chemotaxis model is formulated as follow:
\begin{equation}\label{problem}
\begin{cases}
u_t= \Delta u - \chi \nabla \cdot (u \nabla v)+\xi \nabla \cdot (u \nabla w)  & \textrm{ in } \Omega \times (0,T_{max}),\\
 v_t=\Delta v-\beta v +f(u)  & \textrm{ in } \Omega \times (0,T_{max}),\\
 w_t= \Delta w - \delta w + g(u)& \textrm{ in } \Omega \times (0,T_{max}),\\
u_{\nu}=v_{\nu}=w_{\nu}=0 & \textrm{ on } \partial \Omega \times (0,T_{max}),\\
u(x,0)=u_0(x), \;  v(x,0)= v_0(x), \;  w(x,0)= w_0(x) & x \in \bar\Omega.
\end{cases}
\end{equation}
Herein, $\Omega$ of $\R^n$, with $n\geq 2$, is a bounded and smooth domain, $\chi, \xi, \beta,\delta>0$, and $f=f(s)$ and $g=g(s)$ sufficiently regular functions in their argument $s\geq 0$, essentially behaving as $s^k$ and $s^l$ for some $k,l>0$. Moreover, further regular initial data $u_0(x), v_0(x), w_0(x)\geq 0$ are fixed, $u_\nu$ (and similarly $v_\nu$ and $w_\nu$) indicates the outward normal derivative of $u$ on $\partial \Omega$, whereas $T_{max}$ identifies the maximum time up to which solutions to the system can be extended. 

Once these hypotheses are fixed
\begin{equation}\label{f}
\begin{cases}
f,g \in C^1([0,\infty)) \quad \textrm{with} \quad   0\leq f(s)\leq \alpha s^k  \textrm{ and } \gamma_0 (1+s)^l\leq g(s)\leq \gamma_1 (1+s)^l,\quad  \textrm{for some}\; \alpha, k, l >0, \gamma_1\geq\gamma_0>0,\\
(u_0, v_0,  w_0)\in (W^{1,\infty}(\Omega))^3, \textrm{ with } u_0, v_0,w_0\geq 0 \textrm{ on } \bar{\Omega},
\end{cases}
\end{equation}
\cite[Theorem 2.2]{CFVApplAnal2024} establishes that problem \eqref{problem} admits a unique global and uniformly bounded classical solution (i.e., $\TM=\infty$ and there exists $C>0$ such that $\lVert u(\cdot,t) \rVert_{L^\infty(\Omega)} \leq C$ for all $t\in (0,\infty)$) whenever
\begin{enumerate}[label=(\roman*)]
\item $l,k\in \left(0,\frac{1}{n}\right]$;
\item  $l\in \left(\frac{1}{n},\frac{1}{n}+\frac{2}{n^2+4}\right)$ and $k\in \left(0,\frac{1}{n}\right]$, or  $k\in \left(\frac{1}{n},\frac{1}{n}+\frac{2}{n^2+4}\right)$ and $l\in \left(0,\frac{1}{n}\right]$;
\item  $l,k\in \left(\frac{1}{n},\frac{1}{n}+\frac{2}{n^2+4}\right).$
\end{enumerate} 
From the one hand, we mention that the above conditions have been improved in the recent paper \cite{ColumbuDiazFrassu}; in the specific, \cite[Theorem 2.2]{ColumbuDiazFrassu} ensures boundedness under the more relaxed assumption $k,l\in \left(0,\frac{2}{n}\right).$ 

As to our contribution, we aim at providing a further scenario toward boundedness involving also coefficients connected to $g$ in \eqref{f}; essentially we will show that 
$$\textit{solutions to model \eqref{problem} are uniformly bounded in time whenever $k<l$ and under a largeness assumption on $\gamma_0$}.$$
\begin{remark}[On the origins and the meaning of model \eqref{problem}] 
The interested reader can find motivations connected to biological phenomena described by system \eqref{problem} exactly in \cite{CFVApplAnal2024}, and references therein mentioned. Also known results in close contexts are collected. 
\end{remark}
\section{Presentation of the main theorem}
As an essential tool in order to mathematically formulate our main theorem,  we have first to recall the following consequence of Maximal Sobolev regularity results (\cite{hieber_pruess} or \cite[Theorem 2.3]{giga_sohr}):
\begin{proposition}\label{PropositionConstantC}
For $n\in\N$, let $\Omega\subset \R^n$ be a bounded domain with smooth boundary, $\rho>0$ and $q>\max\{1,\frac{1}{\rho}\}$. Then there is $C_\rho=C_\rho(\Omega,n,q)>0$ such that the following holds: Whenever $T\in(0,\infty]$, $I=[0,T)$, $h\in L^q(I;L^q(\Omega))$ and $\psi_0\in W^{2,q}_{\frac{\partial}{\partial \nu}}(\Omega)= \{\psi_0\in W^{2,q}(\Omega)\,:\, \partial_\nu \psi_0=0 \,\textrm{ on }\, \partial\Omega$\}, every solution $\psi\in W_{loc}^{1,q}(I;L^q(\Omega))\cap L^q_{loc}(I;W^{2,q}(\Omega))$ of 
 \[
  \psi_t=\Delta \psi-\rho \psi + h\;\; \text{ in }\;\;\Omega\times(0,T);\quad
  \partial_{\nu} \psi=0\;\; \text{ on }\;\;\partial\Omega \times(0,T); \quad \psi(\cdot,0)=\psi_0 \;\; \text{ on }\;\;\Omega  
 \]
 satisfies 
 \[
  \int_0^t  e^s{\int_{\Omega} \left(|\psi(\cdot,s)|^q+|\psi_t(\cdot,s)+\frac{\psi(\cdot,s)}{q}|^q+|\Delta \psi(\cdot,s)|^q\right)}ds \le 2^{q-1}C_\rho^q \left[\lVert \psi_0 \rVert_{q,1-\frac{1}{q}}^q+\int_0^t e^s {\int_{\Omega} |h(\cdot,s)|^{q}}ds\right] \quad \text{for all } t\in(0,T).
 \]
 \begin{proof}
The proof is based on the classical result in \cite{PrussSchanaubeltMaximalRegul}; for an appropriate adaptation to our case see details, for instance, in \cite{SachikoEtAlDCDS-B2024}. 
 \end{proof}
\end{proposition}
\begin{remark}[On the constant $C_\rho$ and the norm $\lVert \psi_0 \rVert_{q,1-\frac{1}{q}}$ in Proposition \ref{PropositionConstantC}]
The key role of Proposition \ref{PropositionConstantC} is the existence of the constant $C_\rho$, which  remains defined once $n,\Omega$ and $q$ are set. In particular (see \cite[Theorem 2.5]{PrussSchanaubeltMaximalRegul}), $C_\rho$ does not depend on the initial configuration $\psi_0$ and the source $h$.

As to $\lVert \psi_0 \rVert_{q,1-\frac{1}{q}}$, it represents the norm of $\psi_0$ in the interpolation space $(L^q(\Omega),W^{2,q}_{\frac{\partial}{\partial \nu}(\Omega)})_{1-\frac{1}{q},q}.$ (See, for instance, \cite[$\S$1]{lunardi2012analytic}.)  
\end{remark}
Exactly in view of what said, we can now give the claim of our
\begin{theorem}\label{TheoMain1}
For $n\in \N$, let $\Omega$ be a bounded domain of $\R^n$ with smooth boundary, $0<k<l$, $\delta,\alpha,\beta>0$ and 
\begin{equation}\label{pBar} 
\bar{p}=\max\left\{\frac{n}{2},k\left(\frac{1}{\beta}-1\right),l\left(\frac{1}{\delta}-1\right)\right\}+1.
\end{equation}
Additionally, let us set 
\begin{equation}\label{TehoPositionA-B} 
\mathcal{A}= 2^{-\frac{l(\bar{p}+l-1)+\bar{p}}{\bar{p}+l}}\left(\frac{\bar{p}+l}{\bar{p}+2l+\delta(\bar{p}+l)}\right).
\end{equation}
Then there exists $\mathcal{C}=\mathcal{C}(n,\Omega, l,k,\delta,\beta)>0$ such that if $\mathcal{C}<\mathcal{A}$, it is possible to find $\gamma_1,\gamma_0>0$ fulfilling
\begin{equation}\label{TehoPositionGam_1-Gam_0}
\gamma_1\geq \gamma_0>  \mathcal{A}^{-1}\mathcal{C}\gamma_1
\end{equation}
and  with this  property: Whenever $f,g$, $u_0,v_0,w_0$ are taken as in \eqref{f}, problem \eqref{problem} admits a global and uniformly bounded solution $(u,v,w)\in (C^0(\bar{\Omega}\times [0,\infty))\cap  C^{2,1}(\bar{\Omega}\times (0,\infty)))^3.$ 
\end{theorem}
\section{Existence of local-in-time solutions and a pricniple for boundedness}\label{LocalSol}
The arguments concerning the forthcoming local existence issue and the boundedness criterion are standard; details are achievable in \cite{TaoWanM3ASAttrRep} and \cite[Appendix A.]{TaoWinkParaPara}. 
\subsection{Local existence statement} Once $\chi,\xi,\beta,\delta>0$ and $f,g$,  $u_0, v_0$ are fixed as in \eqref{f}, from here henceforth, with $(u, v, w)$ we will refer to the classical and nonnegative solution to problem \eqref{problem}; $u,v,w$ are defined for all $(x,t) \in \bar{\Omega}\times [0,T_{max})$, for some finite $T_{max}$.
\subsection{Boundedness criterion}
As explained in the next lines, if we establish that $u\in L^\infty((0,\TM);L^p(\Omega))$, for some $p>\frac{n}{2}$,  we can exploit the boundedness criterion below and directly obtain that, indeed,  $u\in L^\infty((0,\infty);L^\infty(\Omega))$; as an immediate consequence of that, well-known parabolic regularity results applied to the equations of $v$ and $w$ entail that also $v,w$ belong to $L^\infty((0,\infty);L^\infty(\Omega)).$ 
 
Definitely, globality and boundedness of $(u,v,w)$, in  the sense that
\begin{align*}
	\begin{split}
		u, v, w\in C^0(\bar{\Omega}\times [0,\infty))\cap  C^{2,1}(\bar{\Omega}\times (0,\infty)) 
		\cap L^\infty((0, \infty);L^{\infty}(\Omega))
	\end{split}
\end{align*}
are achieved whenever this \textit{boundedness criterion} applies:
\begin{equation}\label{BoundednessCriterio}
	\begin{array}{c}
		\textrm{If } \exists  \; L>0, p>\frac{n}{2}\; \Big| \; \displaystyle \int_\Omega u^p \leq L  \textrm{ on }  (0,\TM)\Rightarrow (u,v,w) \in (L^\infty((0, \infty);L^{\infty}(\Omega)))^3.
	\end{array}
\end{equation}
Subsequently, Theorem \ref{TheoMain1} is established once \eqref{BoundednessCriterio} is derived. 
\section{A priori bounds; proof of the main result}
\textit{From now on we will tacitly assume that all the appearing constants below $c_i$, $i=1,2,\ldots,$ are positive.} 
\subsection{Some preparatory tools}\label{PreliminariesSection}
Let us start with this necessary result: 
\begin{lemma}\label{BoundsInequalityLemmaTecnicoYoung}  
Let $A,B \geq 0$ and $p\geq 1$. Then 
we have 
\begin{equation}\label{InequalityA+BToPowerP}
(A+B)^p \leq 2^{p-1}(A^p+B^p).
\end{equation}
\begin{proof}
The proof is available \cite[Theorem 1]{Jameson_2014Inequality}.
\end{proof}
\end{lemma}
\subsection{Achieving the boundedness criterion}
We have this sequence of results, valid for any constant $\Xi>0$, which will be properly chosen later on, in the proof of our theorem.

The following lemma is valid for a general class of proper functions. Despite that, we contextualize it to the local solution $(u,v,w)$ to problem \eqref{problem}.
\begin{lemma}\label{Estim_general_For_u^pLemmaEllittico}
For any $p>1$ and all $t\in (0,\TM)$, we have
 \begin{equation*}
 \begin{split}
 (p-1)\xi \int_\Omega u^p|w_t|& \leq (p-1)\xi \Xi\int_\Omega |w_t+\frac{l}{p+l}w|^\frac{p+l}{l}\\ &
 \quad +\left[\xi p\frac{p-1}{p+l}\left(\Xi\frac{p+l}{l}\right)^{-\frac{l}{p}}\left(1+\frac{l}{p+l}\right)\right]\int_\Omega u^{p+l}+l \xi \Xi \frac{p-1}{p+l}\int_\Omega w^\frac{p+l}{l}, 
 \end{split}
 \end{equation*}
 and
  \begin{equation*}
 (p-1)\xi \delta \int_\Omega u^pw \leq (p-1)\xi \delta\Xi\int_\Omega w^\frac{p+l}{l}+\xi p \delta\frac{p-1}{p+l}\left(\Xi\frac{p+l}{l}\right)^{-\frac{l}{p}}\int_\Omega u^{p+l}. 
 \end{equation*}
\begin{proof}
Let, for commodity but also for reasons which will be clearer later,  $q=\frac{p+l}{l}$. From the evident relation $|w_t|\leq |w_t+\frac{1}{q}w|+|\frac{1}{q}w|,$ we obtain
\begin{equation*}
(p-1)\xi \int_\Omega u^p|w_t|\leq (p-1)\xi \int_\Omega u^p|w_t+\frac{w}{q}|+\xi \frac{p-1}{q}\int_\Omega u^pw \quad \textrm{on } \; (0,\TM),
\end{equation*}
so that thanks to the Young inequality for all $t\in(0,\TM)$ it is seen
\begin{equation*}
\begin{split}
(p-1)\xi \int_\Omega u^p|w_t|&\leq (p-1)\xi \Xi \int_\Omega |w_t+\frac{w}{q}|^q \\ & 
\quad +\frac{(p-1)p \xi}{p+l}(\Xi q)^{-\frac{l}{p}}\int_\Omega u^{p+l}+\Xi \frac{p-1}{q}\xi\int_\Omega w^q +\xi p \frac{p-1}{q(p+l)} \left(\Xi q\right)^{-\frac{l}{p}}\int_\Omega u^{p+l},
\end{split}
\end{equation*} 
and the first claim is established. 

As to the other relation, it can be derived in the same flavor. 
 \end{proof}
 \end{lemma}
 \begin{lemma}\label{Lemma-1}
For any $p>\max\left\{1,l\left(\frac{1}{\delta}-1\right)\right\}$ and $t \in(0,\TM)$ it holds that
\begin{equation}\label{Inequality1ThanksPRegularity} 
\begin{split}
  (p-1)\xi \Xi &\int_0^t e^s\left(\int_\Omega  |w_t+\frac{l}{p+l}w|^\frac{p+l}{l}\right) ds  \leq (p-1)\xi \Xi 2^{\frac{p}{l}}\mathcal{C}_\delta^\frac{p+l}{l} \\ & \times \left[\lVert w_0 \rVert_{\frac{p+l}{l},\frac{p}{p+l}}^\frac{p+l}{l}+\gamma_1^\frac{p+l}{l}2^{p+l-1}\int_0^te^s\left(\int_\Omega u^{p+l}ds\right)+\gamma_1^\frac{p+l}{l}2^{p+l-1}|\Omega|\int_0^t e^s ds\right],  
  \end{split}
 \end{equation}
 and
\begin{equation}\label{Inequality2ThanksPRegularity}
\begin{split}
(p-1)\xi \Xi & \left(\frac{l}{p+l}+\delta\right)\int_0^t e^s\left(\int_\Omega w^\frac{p+l}{l}\right)\leq (p-1)\xi \Xi 2^{\frac{p}{l}} \left(\frac{l}{p+l}+\delta\right)\mathcal{C}_\delta^\frac{p+l}{l}\\ & \times \left[
\lVert w_0 \rVert_{\frac{p+l}{l},\frac{p}{p+l}}^\frac{p+l}{l}+\gamma_1^\frac{p+l}{l}2^{p+l-1}\int_0^te^s\left(\int_\Omega u^{p+l}ds\right)+\gamma_1^\frac{p+l}{l}2^{p+l-1}|\Omega|\int_0^t e^s ds\right].
\end{split}
\end{equation} 
 \begin{proof}
We can derive \eqref{Inequality1ThanksPRegularity} (and similarly \eqref{Inequality2ThanksPRegularity}) by invoking Proposition \ref{PropositionConstantC} with $\psi=w$, $h=g$ and $\rho=\delta$; indeed, for $q=\frac{p+l}{l}$ as before, it is $q>\max\{1,\frac{1}{\delta}\}$ so that 
 \begin{equation*}
  (p-1)\xi \Xi\int_0^t e^s\left(\int_\Omega |w_t+\frac{w}{q}|^q\right) ds \leq (p-1)\xi \Xi \mathcal{C}_\delta^q 2^{q-1}\left[\lVert w_0 \rVert_{q,1-\frac{1}{q}}^q+\int_0^te^s\left(\int_\Omega g(u)^q\right)ds\right],  
 \end{equation*} 
and the conclusion is attained by virtue of the upper bound \eqref{f} for $g$ and  \eqref{InequalityA+BToPowerP}, in the form $(u+1)^{p+l}\leq 2^{p+l-1}(u^{p+l}+1)$ (naturally $p+l>1$).
 \end{proof}
 \end{lemma}
 The next two results, indeed, provide properties of local solutions $(u,v,w)$ to model \eqref{problem} and are based on applications of Proposition \ref{PropositionConstantC}.
 \begin{lemma}\label{Inequality3ThanksPRegularityLemma}
For any $p>\max\left\{1,k\left(\frac{1}{\beta}-1\right)\right\}$ and $t \in(0,\TM)$ there is $\const{321}$ such that
\begin{equation*}
\begin{split}
 \const{321}&\int_0^t e^s\left(\int_\Omega  |\Delta v|^\frac{p+k}{k}\right) ds  \leq \const{321} \mathcal{C}_\beta^{\frac{p+k}{k}} \left[\lVert v_0 \rVert_{\frac{p+k}{k},\frac{p}{p+k}}^\frac{p+k}{k}+\alpha^\frac{p+k}{k}\int_0^te^s\left(\int_\Omega u^{p+k}\right)ds\right].  
  \end{split}
 \end{equation*}
 \begin{proof}
 The proof follows from analogous arguments used in Lemma \ref{Lemma-1}; in this case, in particular, Proposition \ref{PropositionConstantC} is exploited with $q=\frac{p+k}{k}>\max\{1,\frac{1}{\beta}\}$, $\psi=v, h=f$ and $\rho=\beta$.
 \end{proof}
 \end{lemma} 
 With the aim of ensuring $u\in L^\infty((0,\TM);L^p(\Omega))$ for some $p>\frac{n}{2}$, let us study the evolution in time of $t\mapsto\int_\Omega u^p$; this will be done by means of testing procedures. 
\begin{lemma}
 For any $p>1$ and all $t\in (0,\TM)$ the following relation is satisfied:
 \begin{equation}\label{InequalityMxReg4}
 \begin{split}
 \frac{d}{dt}\int_\Omega u^p & \leq  \int_\Omega u^{p+k}+\const{321}\int_\Omega |\Delta v|^\frac{p+k}{k} +(p-1)\xi \Xi \int_\Omega |w_t+\frac{l}{p+l}w|^\frac{p+l}{l}\\&
 +(p-1)\xi \Xi \left(\delta+\frac{l}{p+l}\right)\int_\Omega w^\frac{p+l}{l}+(p-1)\xi \left[\frac{p}{p+l}\left(\frac{p+l}{l}\right)^{-\frac{l}{p}}\Xi^{-\frac{l}{p}}\left(1+\delta+\frac{l}{p+l}\right)-\gamma_0\right]\int_\Omega u^{p+l}.
 \end{split}
 \end{equation}
 \begin{proof}
By testing the first equation of problem \eqref{problem} with $p u^{p-1}$, using its boundary conditions and taking into account the second and the third equation,  we have thanks to the Young inequality
\begin{equation*}
\begin{split}
\frac{d}{dt} \int_\Omega u^p &=p \int_\Omega u^{p-1}u_t
= -p (p-1) \int_\Omega u^{p-2} |\nabla u|^2 - (p-1) \chi \int_\Omega u^p \Delta v+ (p-1)\xi \int_\Omega u^p \Delta w\\
&\leq  \int_\Omega u^{p+k}+\const{321}\int_\Omega |\Delta v|^\frac{p+k}{k} +(p-1)\xi \int_\Omega u^p (w_t+\delta w -g(u))\quad \textrm{on }\, (0,\TM).
\end{split}  
\end{equation*}  
Now, we recall the properties of  $g$ given in \eqref{f} so to deduce, by using Young's inequality, again the relation $|w_t|\leq |w_t+\frac{l}{p+l}w|+|\frac{l}{p+l}w|,$ and Lemma \ref{Estim_general_For_u^pLemmaEllittico}
\begin{equation*}
\begin{split}
\frac{d}{dt} \int_\Omega u^p &\leq  \int_\Omega u^{p+k}+\const{321}\int_\Omega |\Delta v|^\frac{p+k}{k} +(p-1)\xi \int_\Omega u^p |w_t|+\delta (p-1)\xi \int_\Omega u^p w  -\xi \gamma_0 (p-1) \int_\Omega u^{p+l}\\&
\leq  \int_\Omega u^{p+k}+\const{321}\int_\Omega |\Delta v|^\frac{p+k}{k}+(p-1)\xi \Xi \int_\Omega |w_t+\frac{l}{p+l}w|^\frac{p+l}{l}\\&
\quad +\frac{l(p-1)\xi \Xi}{p+l}\int_\Omega w^\frac{p+l}{l}
+\frac{p(p-1)\xi}{p+l}\left(\Xi \frac{p+l}{l}\right)^{-\frac{l}{p}}\left(1+\frac{l}{p+l}\right)\int_\Omega u^{p+l}\\ &
\quad +(p-1)\xi \delta \Xi\int_\Omega w^\frac{p+l}{l}
+\frac{p(p-1)\xi \delta }{p+l}\left(\Xi \frac{p+l}{l}\right)^{-\frac{l}{p}}\int_\Omega u^{p+l}-\xi (p-1)\gamma_0 \int_\Omega u^{p+l}
\quad \textrm{on }\, (0,\TM).
\end{split}
\end{equation*}
The claim is achieved by collecting terms. 
\end{proof}
\end{lemma}
\begin{lemma}\label{Estim_general_For_u^pLemmaParabolico2} 
Let $k<l$. Then for every $p>\max\left\{\frac{n}{2},k\left(\frac{1}{\beta}-1\right),l\left(\frac{1}{\delta}-1\right)\right\}$  we have that for all $t<\TM$
 \begin{equation*}
 \begin{split}
e^t\int_\Omega u^p  \leq & \const{ga}+ \const{ga1} \int_0^t e^sds + \\&
 (p-1)\xi \left[\left(1+\delta+\frac{l}{p+l}\right)\left( \Xi C_\delta^\frac{p+l}{l} \gamma_1^{\frac{p+l}{l}}2^{\frac{p}{l}+p+l-1}+\frac{p}{p+l}\left(\frac{p+l}{l}\right)^{-\frac{l}{p}}\Xi^{-\frac{l}{p}}\right)+\varepsilon-\gamma_0\right]\int_0^t e^s\left(\int_\Omega u^{p+l}\right)ds.
 \end{split}
 \end{equation*}  
\begin{proof}
Let us start with these estimates, fruit of the application of Young's inequality: for all $\varepsilon>0$, $\hat{c}>0$, $p>1$ and $0<k<l$ it holds that 
\begin{equation}\label{u^pTou^p+l}
\hat{c} \int_\Omega u^p\leq \frac{\varepsilon}{2} \int_\Omega u^{p+l}+\const{uuu} \quad \textrm{for all } t\in (0,\TM),
\end{equation}  
and
\begin{equation}\label{u^p+kTou^p+l}
\hat{c} \int_\Omega u^{p+k}\leq \frac{\varepsilon}{2} \int_\Omega u^{p+l}+\const{uut} \quad \textrm{for all } t\in (0,\TM).
\end{equation}  
By adding to both sides of relation  \eqref{InequalityMxReg4} the term $\int_\Omega u^p$, estimate \eqref{u^pTou^p+l} leads to this inequality, valid on $(0,\TM)$.
 \begin{equation}\label{AAA}
 \begin{split}
 \frac{d}{dt}\int_\Omega u^p&  +\int_\Omega u^p  \leq   \int_\Omega u^{p+k}+\const{321}\int_\Omega |\Delta v|^\frac{p+k}{k} +(p-1)\xi \Xi \int_\Omega |w_t+\frac{l}{p+l}w|^\frac{p+l}{l}\\&
 +(p-1)\xi \Xi \left(\delta+\frac{l}{p+l}\right)\int_\Omega w^\frac{p+l}{l}+(p-1)\xi \left[\frac{p}{p+l}\left(\frac{p+l}{l}\right)^{-\frac{l}{p}}\Xi^{-\frac{l}{p}}\left(1+\delta+\frac{l}{p+l}\right)+\frac{\varepsilon}{2}-\gamma_0\right]\int_\Omega u^{p+l}+\const{xxx}.
 \end{split}
 \end{equation}
 Successively, we multiply \eqref{AAA} by $e^t$ and integrate on $(0,t)$. From the identity $\frac{d}{dt}(e^t\int_\Omega u^p) =e^t\frac{d}{dt}\int_\Omega u^p +e^t\int_\Omega u^p$, we get
 \begin{equation*}
 \begin{split}
e^t&\int_\Omega u^p  \leq \int_\Omega u_0^p+   \int_0^t e^s\left\{\int_\Omega u^{p+k}+\const{321}\int_\Omega |\Delta v|^\frac{p+k}{k} +(p-1)\xi \Xi \int_\Omega |w_t+\frac{l}{p+l}w|^\frac{p+l}{l}\right.\\&
\left. +(p-1)\xi \Xi \left(\delta+\frac{l}{p+l}\right)\int_\Omega w^\frac{p+l}{l}+(p-1)\xi\left[\frac{p}{p+l}\left(\frac{p+l}{l}\right)^{-\frac{l}{p}} \Xi^{-\frac{l}{p}}\left(1+\delta+\frac{l}{p+l}\right)+\frac{\varepsilon}{2}-\gamma_0\right]\int_\Omega u^{p+l}+\const{xxx}\right\}ds.
 \end{split}
 \end{equation*} 
 The term involving $\int_0^t e^s(\int_\Omega |\Delta v|^{\frac{p+k}{k}})ds$  can be essentially controlled  by $\int_0^t e^s(\int_\Omega u^{p+k})ds$, thanks to Lemma \ref{Inequality3ThanksPRegularityLemma}; additionally, $\int_\Omega u^{p+k}$ is treated through \eqref{u^p+kTou^p+l}. These two operations provide
   \begin{equation*}
 \begin{split}
e^t\int_\Omega u^p  \leq & \int_\Omega u_0^p+   \int_0^t e^s\left\{(p-1)\xi \Xi \int_\Omega |w_t+\frac{l}{p+l}w|^\frac{p+l}{l}+(p-1)\xi \Xi \left(\delta+\frac{l}{p+l}\right)\int_\Omega w^\frac{p+l}{l}\right.\\&
\left.+ (p-1)\xi \left[\frac{p}{p+l}\left(\frac{p+l}{l}\right)^{-\frac{l}{p}}\Xi^{-\frac{l}{p}}\left(1+\delta+\frac{l}{p+l}\right)+\varepsilon-\gamma_0\right]\int_\Omega u^{p+l}+\const{xxxy}\right\}ds.
 \end{split}
 \end{equation*} 
 By invoking \eqref{Inequality1ThanksPRegularity} and \eqref{Inequality2ThanksPRegularity} we conclude by virtue of a reorganization of the involved terms.  
 \end{proof}
 \end{lemma}
 \subsection{Proof of Theorem \ref{TheoMain1}} 
With the above preparations we are now in a position to establish what anticipated.
\begin{proof}
 For $0<k<l$, $\delta,\alpha,\beta>0$, let $\bar{p}$ be as in \eqref{pBar} and, additionally, for $C_\delta\left(\Omega,n,\frac{\bar{p}+l}{l}\right)$ being the constant provided by Proposition \ref{PropositionConstantC}, when it is applied to the equation for $w$ in model \eqref{problem}, let also set
%
\begin{equation*}
\mathcal{C}=C_\delta\left(\Omega,n,\frac{\bar{p}+l}{l}\right).
\end{equation*}
Since by assumptions $\mathcal{C}<\mathcal{A}$, where $\mathcal{A}$ is defined in \eqref{TehoPositionA-B}, we can find $\gamma_1\geq \gamma_0$ complying with  
\eqref{TehoPositionGam_1-Gam_0}; in these positions, let
\begin{equation*}
\Xi= \frac{l}{\bar{p}+l}\mathcal{C}_\delta^{-\frac{\bar{p}}{l}}\gamma_1^{-\frac{\bar{p}}{l}}2^{-\frac{\bar{p}(\bar{p}+(\bar{p}+l-1)l)}{l(\bar{p}+l)}},
\end{equation*}
and let $f,g,u_0,v_0$ and $w_0$ obey \eqref{f}. Some computations show that for proper small $\varepsilon>0$ 
\begin{equation*}
 \gamma_0>\mathcal{A}^{-1}\mathcal{C}\gamma_1 \Rightarrow \left(1+\delta+\frac{l}{\bar{p}+l}\right)\left( \Xi C_\delta^\frac{\bar{p}+l}{l} \gamma_1^{\frac{\bar{p}+l}{l}}2^{\frac{\bar{p}}{l}+\bar{p}+l-1}+\frac{\bar{p}}{\bar{p}+l}\left(\frac{\bar{p}+l}{l}\right)^{-\frac{l}{\bar{p}}}\Xi^{-\frac{l}{\bar{p}}}\right)+\varepsilon-\gamma_0\leq 0,
\end{equation*}
and henceforth hypothesis \eqref{TehoPositionGam_1-Gam_0} allows to exploit Lemma \ref{Estim_general_For_u^pLemmaParabolico2} and obtain 
 \begin{equation*}
e^t\int_\Omega u^{\bar{p}}  \leq  \const{ga}+ \const{ga1} \int_0^t e^sds \quad \textrm{for all } t\in (0,\TM),
 \end{equation*}
or also $\int_\Omega u^{\bar{p}}\leq L$ on $(0,\TM)$. The claim follows from the extensibility criterion \eqref{BoundednessCriterio}.
\end{proof}

\subsubsection*{Acknowledgments}
The authors are members of the Gruppo Nazionale per l'Analisi Matematica, la Probabilit\`a e le loro Applicazioni (GNAMPA) of the Istituto Na\-zio\-na\-le di Alta Matematica (INdAM). The authors are partially supported by the research projects \textit{Analysis of PDEs in connection with
real phenomena}, funded by Fondazione di Sardegna (2021, CUP F73C22001130007) and   MIUR (Italian Ministry of Education,
University and Research) Prin 2022 \textit{Nonlinear differential problems with applications to real phenomena} (Grant Number:
2022ZXZTN2). SF also acknowledges financial support by INdAM-GNAMPA project \textit{Problemi non lineari di tipo stazionario ed evolutivo} (CUP E53C23001670001).

\end{document}